\documentclass[12pt]{article}
\usepackage{amssymb}
\usepackage{array}
\usepackage{latexsym}
\usepackage{enumerate}
\usepackage{amsmath}
\usepackage{amsfonts}
\usepackage{amsthm}
\usepackage{psfrag}
\usepackage{comment}
\usepackage[english]{babel}

\title{\bf{Nakajima's remark on Henn's proof}}
\date{}
\author{M.~Giulietti ${}^*$ and G.~Korchm\'aros
 \thanks{Research supported by  the Italian
    Ministry MURST, Strutture geometriche, combinatoria e loro
    applicazioni.}}

{\theoremstyle{definition}
\newtheorem*{definition*}{Definition}

\newtheorem*{proposition*}{Proposition}
\newtheorem*{corollary*}{Corollary}
\newtheorem*{lemma*}{Lemma}

\def\Aut{\mbox{\rm Aut}}

\def\deg{\mbox{\rm deg}}

\def\Aut{\mbox{\rm Aut}}

\def\Div{\mbox{\rm div}}

\newcommand{\ga}{\alpha}

\newfont{\german}{eufm10 scaled 1200}
\newfont{\sgerman}{eufm10 scaled 900}

\newcommand{\GB}{\mbox{\german B}}

\sloppy

\begin{document}
\maketitle

    \begin{abstract} We fill up a gap in Henn's proof concerning large
    automorphism groups of function fields of degree $1$ over an
    algebraically closed field of positive characteristic.
    \end{abstract}
    \section{Introduction} In 1973, Stichtenoth \cite{stichtenoth1973I}
showed that the Hermitian function fields are the unique function
fields $K$ of transcendency degree $1$ over an algebraically
closed ground field $\Omega$ of characteristic $p$ whose
automorphism group $G=\Aut(K/\Omega)$ has order at least $16g^4$
where $g\geq 2$ is the genus of $K$.

In 1978, Henn \cite{henn1978} gave a complete classification of
such function fields under the weaker hypothesis $|G|\geq 8g^3$.
Later Nakajima \cite{nakajima1987} improved Henn's result for
ordinary curves. In a footnote of his paper, Nakajima claimed:
``Stichtenoth's result was improved by Henn. But his proof
contains a gap (last paragraph of pg.\,104). I do not know if the
gap can be covered.'' The gap appears to be in the penultimate
line on pg.\, 104 when Henn claims ``und ersichtlich $z\cong
\zeta(z)$ gilt, folgt hieraus $E\le 2$''. The purpose of the
present note is to fill up this gap. Actually, we are going to
show the missing details in Henn's proof.

We keep notation and terminology from \cite{henn1978}. In the last
paragraph on pg.\,104, the following case is investigated: $K$ has
two places $\GB$ and $\GB_2$ with $G_0(\GB)=G_0(\GB_2)$ but none
of the hypotheses in Proposition 1 holds.

For this case, Henn explicitly proves the following claims:
\begin{enumerate}
    \item[(i)] $\nu_1=2$;
    \item[(ii)] the genus $g_2$ of $K^{G_2(\GB)}$ is equal to
    $q_1-1$ with $q_1=|G_1(\GB)/G_2(\GB)|$.
    \item[(iii)] the smallest non-gap at $\widetilde{\GB}$ is $q_1$.
\end{enumerate}

Then he observes that $q_1+1$ must also be a pole number at
$\widetilde{\GB}$, and chooses an element $z\in K$ such that
$$z\cong \frac{\widetilde{\GB}_2^a {\textsf{a}}}{\widetilde{\GB}^{q_1+1}},\quad
a>1,\,\, \widetilde{\GB}_2 \nmid {\textsf{a}}.$$ For $\sigma\in
{\rm{Gal}}(K^{G_2(\GB)}/K^{G_1(\GB)})$, he shows that
$\sigma(z)=z+\alpha$ with $\alpha\in \Omega$ and $\alpha\neq 0$.
From this and Lemma 2, he deduces that
$$2=\nu_1=1-a+q_1, $$
whence $a=q_1-1$ and $\deg ({\textsf{a}})=2$ follow. At this point
Henn takes an element $\zeta$ of order $|\zeta|=E$ from $G_0(\GB)$.
To end the proof it is sufficient to show that $\zeta$ has order at
most $2$, as the hypothesis $1<e<E$ will give then a contradiction.
Henn's idea is to show first that
\begin{equation}\label{zzz} \zeta(z)=cz\,\,\text{with } c\in
\Omega\setminus \{0\},
\end{equation}
 and then to deduce $E=2$ from it. He claims
that (\ref{zzz}) follows from the fact that ``jeder Punkt $\neq \GB,
\GB_2$ unter $\zeta$ genau $E$ Konjugierte hat.''

We are going to show that $z$ may be chosen such a way that
(\ref{zzz}) holds indeed.

Since $q_1$ and $q_1+1$ are coprime, we have
$K^{G_2(\GB)}=\Omega(y,z)$ where $y$ as in Henn's paper has the
property: $$y\cong
\frac{\widetilde{\GB}_2^{q_1}}{\widetilde{\GB}^{q_1}}.$$

Let $f\in \Omega[Z,Y]$ be an irreducible polynomial over $\Omega$
such that $f(z,y)=0$. It may be noted that the plane irreducible
curve $C$ with affine equation $f(Z,Y)=0$ is in its Weierstrass
normal form:
$$f(Z,Y)=Y^{q_1+1}+\gamma\,Z^{q_1}+U_1(Z)Y^{q_1}+\ldots+U_{q_1}(Z)Y+U_{q_1+1}(Z),$$
where $\gamma\in \Omega\setminus \{0\}$ and $\deg\,U_i(Z)\leq
iq_1/(q_1+1)$ for $i=1,\ldots, q_1$ and $\deg\,U_{q_1+1}(Z)<q_1$.

In particular, $C$ has only one point at infinity, namely the
infinity point $Z_\infty$ of the $Z$-axis which is the center of
the place $\widetilde{\GB}$. Furthermore, the origin $O$ is the
center of $\widetilde{\GB}_2$, and no other place of
$K^{G_2(\GB)}$ is centered either at $Z_\infty$ or $O$.

Since $\zeta$ belongs to the normalizer of $G_2(\GB)$,  $\zeta$
may be viewed as an $\Omega$-automorphism of $K^{G_2(\GB)}$ of
order $E$. Moreover, $\zeta$ is a linear collineation that
preserves $C$. Since $\zeta$ fixes both $\widetilde{\GB}$ and
$\widetilde{\GB}_2$, there exist $\alpha,\beta,\gamma\in \Omega$
with $\ga,\gamma\neq 0$, such that
\begin{equation*}
\begin{array}{ll}
\zeta(y)= & \alpha y; \\
\zeta(z)= & \beta y + \gamma z.
\end{array}
\end{equation*}
If $\beta=0$, then (\ref{zzz}) holds. Now assume that $\beta\neq 0$.
Then $\alpha\neq \gamma$. In fact, if $\alpha=\gamma$ then
$\zeta^j(y)=\alpha^j(y)$ and $\zeta^{j}(z)=j\alpha^{j-1}\beta
y+\alpha^j z$ for every positive integer $j$. Now, for $j=E$ this
implies that $\alpha^E=1$ and hence that $E\equiv 0 \pmod p$, a
contradiction.

Let $u=\beta/(\alpha-\gamma)$, and $z'=z-uy$. Then
$\zeta(z')=\zeta(z)-u\zeta(y)=\gamma z +\beta y - u \alpha y=\gamma
(z-uy)$. Replacing $z$ by $z'$ we obtain
$$\zeta(z')=\gamma z'.$$
  Actually $z$ may be replaced by $z'$ from the
very beginning of the argument, therefore Henn's claim (\ref{zzz})
may be assumed to be true.

Henn's also claims without proof that (\ref{zzz}) implies  $E\le 2$.
This can be shown as follows. From
 $\Div
\zeta(z)=\Div z$ it follows that $\zeta$ must preserve the above
divisor ${\textsf{a}}$. Since $\deg \,{\textsf{a}}=2$, this implies
that $\zeta^2$ fixes each place in the support of ${\textsf{a}}$.
Therefore, $\zeta^2$ as an $\Omega$-automorphism of the rational
function field $K^{G_1(\GB)}$ fixes at least three distinct places,
namely, $\GB',\GB_2'$ and each place lying under those in the
support of ${\textsf{a}}$ in the covering $K^{G_2(\GB)}\to
K^{G_1(\GB)}$. But then $\zeta^2$ is the identity, and hence $E=2$.

\section*{Remark} A revised proof of Henn's classification is found in \cite[Chapter
11.12]{hirschfeld-korchmaros-torres2008}.

\section*{Acknowledgements}
We are grateful to Michel Matignon for turning our attention to
Nakajima's remark about Henn's paper.

    \end{document}